\documentclass[graybox]{svmult}

\usepackage{type1cm}        
\usepackage{makeidx}         
\usepackage{graphicx}        
\usepackage{multicol}        
\usepackage[bottom]{footmisc}
\usepackage{hyperref}
\usepackage{newtxtext}       %
\usepackage{newtxmath}       

\makeindex             

\begin{document}

\title*{On the nonlinear Dirichlet-Neumann method and
    preconditioner for Newton's method}
\author{F. Chaouqui, M. J. Gander, P. M. Kumbhar and T. Vanzan}
\institute{F. Chaouqui \at Temple University, Philadelphia, USA, \email{faycal.chaouqui@temple.edu},
\and M. J. Gander \at Universit\'e de Gen\`eve, Switzerland \email{martin.gander@unige.ch},
\and P. M. Kumbhar  \at Karlsruher Institut f\"{u}r Technologie, Germany \email{pratik.kumbhar@kit.edu},
\and  T. Vanzan \at \'{E}cole Polytecnique F\'{e}d\'{e}rale de Lausanne, Switzerland \email{tommaso.vanzan@epfl.ch}}
%
%
\maketitle

\abstract*{The Dirichlet-Neumann (DN) method has been extensively
  studied for linear partial differential equations, while little
  attention has been devoted to the nonlinear case. In this paper, we
  analyze the DN method both as a nonlinear iterative method and as a
  preconditioner for Newton's method. We discuss the nilpotent
  property and prove that under special conditions, there exists a
  relaxation parameter such that the DN method converges
  quadratically. We further prove that the convergence of Newton's
  method preconditioned by the DN method is independent from the
  relaxation parameter. Our numerical experiments further illustrate
  the mesh independent convergence of the DN method and compare it
  with other standard nonlinear preconditioners.}

\section{Introduction}\label{kumbhar_p_mini_17_sec:intro}

We consider a nonlinear Partial Differential
Equation (PDE)
\begin{equation} \label{kumbhar_p_mini_17_eq:model_gen}
  \mathcal{L}(u)=f\quad \text{in}\quad \Omega,\quad
     u=g\quad \text{on } \partial \Omega,
\end{equation}
where $\Omega \subset \mathbb{R}^d$ for $d\in \left\{1,2,3\right\}$ is
an open bounded domain with a polygonal boundary $\partial \Omega$,
and $f,g \in L^2(\Omega)$.
We suppose that \eqref{kumbhar_p_mini_17_eq:model_gen} admits a unique weak
solution in some Hilbert space $u\in \mathcal{X}$(
e.g. $H^1(\Omega))$. For instance, for a quasilinear operator
  $\mathcal{L}$ in divergence form, explicit assumptions can
  be found in \cite{kumbhar_p_mini_17_cai_dryja_1994} and references therein, see also
  \cite[Chapter 8-9]{kumbhar_p_mini_17_evans2010partial} and \cite[Chapter 9]{kumbhar_p_mini_17_Ciarlet}. 
  Let us divide $\Omega$ into two nonoverlapping subdomains $\Omega_1$
and $\Omega_2$ and define $\Gamma_j=\partial
\Omega_j\setminus\partial\Omega$, $j=1,2$. Let $u_j$ be the
  restriction of $u$ to $\Omega_j$. The nonlinear Dirichlet-Neumann
(DN) method starts from an initial guess $\lambda^0$ and computes for
$n\geq 1$ until convergence
\begin{equation}\label{kumbhar_p_mini_17_eq:DN_method}
\begin{array}{l l l l l}
\mathcal{L}(u_1^n)&=f_1,\quad \text{in}\quad \Omega_1, \quad & &\quad  \mathcal{L}(u_2^n)&=f_2,\quad \text{in}\quad \Omega_2 \\
u^n_1&=g_1, \quad \text{on}\quad \partial \Omega_1\setminus \Gamma, \quad& &\quad u^n_2&=g_2, \quad \text{on}\quad \partial \Omega_2\setminus \Gamma\\
u^n_1&=\lambda^{n} \quad \text{on}\quad  \Gamma,\quad & & \quad \mathcal{N}_2u^n_2&=-\mathcal{N}_1u^n_1 \quad \text{on}\quad  \Gamma,
\end{array}
\end{equation}
where $\lambda^{n}=(1-\theta)\lambda^{n-1}+\theta u^{n-1}_{2
  |\Gamma},$ with $\theta\in (0,1)$ , $f_j:=f|_{\Omega_j}$ 
  and $g_j:=g|_{\partial \Omega_j\setminus \Gamma}$ for $j=1,2$.  The operators $\mathcal{N}_j$
represent the outward nonlinear Neumann conditions that must be imposed
on the interface $\Gamma$ and are usually found through
integration by parts of the variational formulation of the PDE. For
instance, if $\mathcal{L}(u)=-\partial_x ((1+\alpha u^2)\partial_x
u)$, then $\mathcal{N}_ju= (-1)^{j+1} (1+ \alpha
u_{|\Gamma}^2)\partial_x u_{|\Gamma}$. For the well-posedness of the Dirichlet-Neumann method, we further assume that $\mathcal{N}_j u$ defines a bounded linear functional over $\mathcal{X}$.

System \eqref{kumbhar_p_mini_17_eq:DN_method} can be formulated as an iteration over the
substructured variable $\lambda$ as
\begin{equation}\label{kumbhar_p_mini_17_eq:DN_method_SUB}
\lambda^n=G(\lambda^{n-1}):= (1-\theta)\lambda^{n-1} + \theta \text{NtD}_2 \left(-\text{DtN}_1\left(\lambda^{n-1},\psi_1\right),\psi_2\right),
\end{equation}
where $\psi_j=(f_j,g_j)$, $j=1,2$, represent the force term and
boundary conditions, while the nonlinear Dirichlet-to-Neumann
($\text{DtN}_j$) and Neumann-to-Dirichlet operators ($\text{NtD}_j$)
are defined as $\text{DtN}_j(\lambda,\psi_j):=\mathcal{N}_ju_j$, and
$\text{NtD}_j(\phi,\psi_j):=v_j|_\Gamma$, with
\begin{equation}\label{kumbhar_p_mini_17_eq:definition_DtN_NtD}
\begin{array}{rcllrcll}
  \mathcal{L}(u_j)&=&f_j\quad&\text{in $\Omega_j$}, &
  \mathcal{L}(v_j)&=&f_j\quad&\text{in $\Omega_j$}, \\
  u_j&=&g_j \quad &\text{on $\partial \Omega_j\setminus \Gamma$},\quad &
  v_j&=&g_j \quad &\text{on $\partial \Omega_j\setminus \Gamma$},\\
  u_j&=&\lambda& \text{on $\Gamma$} &
  \mathcal{N}_jv_j&=&\phi \quad &\text{on  $\Gamma$}.
\end{array}
\end{equation}
If $u_{\text{ex}}\in H^1(\Omega)$ is the solution of \eqref{kumbhar_p_mini_17_eq:model_gen},
then it must have continuous Dirichlet trace and Neumann flux along
the interface $\Gamma$.  Defining $u_\Gamma:=u_{\text{ex}}|_\Gamma$,
$\phi:=\mathcal{N}_1u_{\text{ex}}|_\Gamma$ and using the operators
$\text{DtN}_j$ and $\text{NtD}_j$, these necessary properties are
equivalent to
\begin{equation}\label{kumbhar_p_mini_17_eq:condition_optimality}
\text{DtN}_1(u_\Gamma,\psi_1)=-\text{DtN}_2(u_\Gamma,\psi_2),\quad\text{and}\quad \text{NtD}_1(\phi,\psi_1)=\text{NtD}_2(-\phi,\psi_2).
\end{equation}

\section{Nilpotent property and quadratic convergence}\label{kumbhar_p_mini_17_sec:Nil_quad}    

It is well known, see e.g. \cite{kumbhar_p_mini_17_quarteroni1999domain,kumbhar_p_mini_17_bookCG}, that if
$\mathcal{L}$ is linear and the subdomain decomposition is symmetric,
then the DN method converges in one iteration for
$\theta=1/2$. Indeed, if $\mathcal{L}$ is linear, one can work on the
error equation, i.e. $\psi_j=0$, and the symmetry of the decomposition
is sufficient to guarantee
$\text{DtN}_1(\cdot,0)\equiv\text{DtN}_2(\cdot,0)$, so that
\begin{equation}\label{kumbhar_p_mini_17_eq:direct_method_linear}
\begin{aligned}
\lambda^1&=\frac{1}{2}\left( \lambda^0 +\text{NtD}_2\left(-\text{DtN}_1(\lambda^0,0)\right),0\right)=\frac{1}{2}\left( \lambda^0 + \text{NtD}_2\left(-\text{DtN}_2(\lambda^0,0)\right),0\right) \\
&=\frac{1}{2}( \lambda^0 - \text{NtD}_2(\text{DtN}_2(\lambda^0,0),0)=0,
\end{aligned}
\end{equation}
where in the third equality we used linearity, and in the last
$\text{NtD}_2(\text{DtN}_2(\lambda,\psi),\psi)=\lambda$.
Can the nonlinear DN method also converge in one iteration?

On the one hand, the relation
$\text{NtD}_j(\text{DtN}_j(\lambda,\psi),\psi)=\lambda$ holds even in
the nonlinear case, simply because the nonlinear $\text{DtN}_{j}$
operator is the inverse of the nonlinear $\text{NtD}_j$ operator.  On
the other hand, due to the nonlinearity of $\mathcal{L}$, one cannot
rely on the error equation, cannot state that
$\text{NtD}_2(-\phi)=-\text{NtD}_2(\phi)$, and the symmetry of the
decomposition is not sufficient to guarantee
$\text{DtN}_1(\lambda,\psi_1)\equiv\text{DtN}_2(\lambda,\psi_2)$,
because of the boundary conditions and the force term.

\begin{figure}[t]
  \centering
  \includegraphics[width=0.48\textwidth,clip]{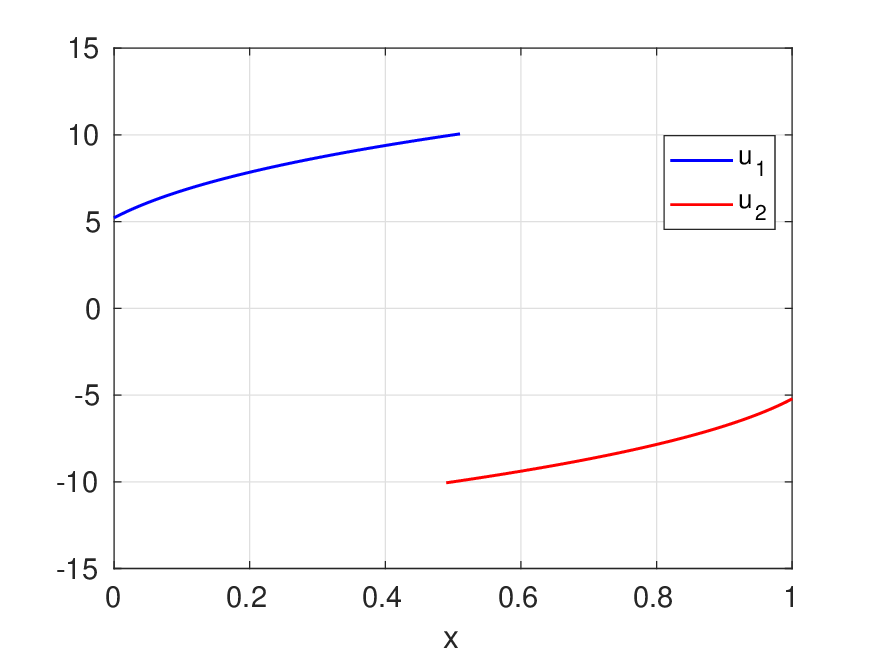}\quad
  \includegraphics[width=0.48\textwidth,clip]{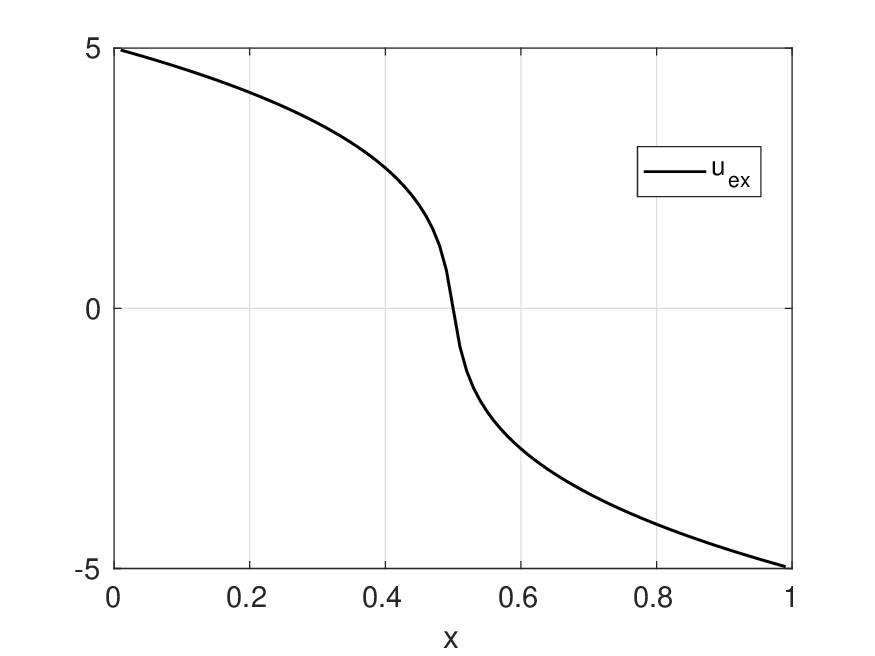}
  \caption{Subdomain solutions of the nonlinear DN method after one
    iteration (left), and exact solution (right). The parameters are
    $g=5$ and $k=2$.}\label{kumbhar_p_mini_17_Fig:nilpotent}
\end{figure}

A straight forward observation is that if the nonlinear DN method
converges in one iteration, then $G(\lambda)=\lambda_{\text{ex}}$, $\forall
\lambda$, that is $G(\cdot)$ is a constant.  A necessary and
sufficient condition for the nonlinear DN method to converge
in one iteration is then
\begin{equation}\label{kumbhar_p_mini_17_eq:conditions_nilpotent}
  0\!=\!G^\prime(\lambda)
  \!=\! \frac{1}{2} +\frac{1}{2}\!\left(\text{NtD}_2
    \left(-\text{DtN}_1(\lambda,\psi_1),\psi_2\right)\right)^\prime
    \!\implies\! \left(\text{NtD}_2\!
    \left(-\text{DtN}_1(\lambda,\psi_1),\psi_2\right)\right)^\prime\!=\!-1.
\end{equation}
Clearly, \eqref{kumbhar_p_mini_17_eq:conditions_nilpotent} is satisfied if
$\text{NtD}_2(-\text{DtN}_1(\lambda,\psi_1),\psi_2)=-\lambda$.  We
consider a toy example in which this condition is satisfied. Let
$\mathcal{L}=-\partial_x\left((1+u^2)\partial_x u\right)$,
$u(0)=g\in\mathbb{R}^+$, $u(1)=-g$ and $f(x)=\sin((2k)\pi x)$. On the
left plot of Fig. \ref{kumbhar_p_mini_17_Fig:nilpotent},
we show the subdomain solutions $u_1$ and $u_2$ obtained from
\eqref{kumbhar_p_mini_17_eq:DN_method} after the first iteration. The two contributions
sum to zero, which is the value of $\lambda_{\text{ex}}$. Thus, after one
iteration we obtain the exact solution shown in the right panel.

Even though the nilpotent property does not hold in general, we show
in the following Theorem that the nonlinear DN method can exhibit
quadratic convergence.
\begin{theorem}[Quadratic convergence of nonlinear DN]\label{kumbhar_p_mini_17_thm:quadratic_convergence}
For any one-dimensional nonlinear problem $\mathcal{L}(u)=f$ such that
$\text{DtN}_1^\prime(\lambda_{\text{ex}},\psi_1)\cdot
\text{DtN}_2^\prime(\lambda_{\text{ex}},\psi_2)>0$ with
$\lambda_{\text{ex}}:=u_{\text{ex}}|_\Gamma$, there exists a $\theta \in (0,1)$ such
that the nonlinear Dirichlet-Neumann method converges quadratically.
\end{theorem}
\begin{proof}
A sufficient condition for quadratic convergence is that the Jacobian of $G(\cdot)$, defined in \eqref{kumbhar_p_mini_17_eq:DN_method_SUB}, is zero at $\lambda_{\text{ex}}:=u_{\text{ex}}|_\Gamma$, that is $G^\prime(\lambda_{\text{ex}})=0$. A direct calculation shows
\begin{equation}\label{kumbhar_p_mini_17_eq:derivative_of_G}
G^\prime(\lambda)= (1-\theta) + \theta \text{NtD}_2^\prime\left(-\text{DtN}_1(\lambda,\psi_1),\psi_2\right)\cdot\left(- \text{DtN}_1^\prime(\lambda,\psi_1)\right).
\end{equation}
Setting $\lambda=\lambda_{\text{ex}}$ and using the optimality condition $\text{DtN}_1(\lambda_{\text{ex}},\psi_1)=-\text{DtN}_2(\lambda_{\text{ex}},\psi_2)$ of \eqref{kumbhar_p_mini_17_eq:condition_optimality},
the above equation changes to
\begin{equation}\label{kumbhar_p_mini_17_eq:derivative_of_G_in_ex}
G^\prime(\lambda_{\text{ex}}) =(1-\theta) + \theta \text{NtD}_2^\prime\left(\text{DtN}_2(\lambda_{\text{ex}},\psi_2),\psi_2\right)\cdot\left(- \text{DtN}_1^\prime(\lambda_{\text{ex}},\psi_1)\right).
\end{equation}
If
$\text{DtN}^\prime_1(\lambda_{\text{ex}},\psi_1)=\text{DtN}^\prime_2(\lambda_{\text{ex}},\psi_2)$
held true, then using the identity
$\text{NtD}_2^\prime\left(\text{DtN}_2(\lambda,\psi_2),\psi_2\right)\cdot\left(\text{DtN}_2^\prime(\lambda,\psi_2)\right)=1$,
obtained by differentiating
$\text{NtD}_j(\text{DtN}_j(\lambda,\psi_j),\psi_j)=\lambda$, we would
easily get that $\theta=1/2$ leads to $G^\prime(\lambda_{\text{ex}})=0$.
Nevertheless, variational calculus shows that to calculate
$\text{DtN}_j^\prime(\lambda,\psi_j)$, one has to solve a linear PDE
which does not depend on $\psi_j$ anymore, but whose coefficients
still depend on the subdomain solutions $u_{\text{ex}}|_{\Omega_1}$ and
$u_{\text{ex}}|_{\Omega_2}$. In general then,
$\text{DtN}^\prime_1(\lambda_{\text{ex}},\psi_1)\neq\text{DtN}^\prime_2(\lambda_{\text{ex}},\psi_2)$.
However, $\text{DtN}_j$ being one dimensional functions, we have
$\text{DtN}^\prime_1(\lambda_{\text{ex}},\psi_1)=\delta
\text{DtN}^\prime_2(\lambda_{\text{ex}},\psi_2)$, for some $\delta\in
\mathbb{R}^+$ if $\text{DtN}_1^\prime(\lambda_{\text{ex}},\psi_1)\cdot
\text{DtN}_2^\prime(\lambda_{\text{ex}},\psi_2)>0$.
 Inserting this into
\eqref{kumbhar_p_mini_17_eq:derivative_of_G_in_ex}, we obtain $G^\prime (\lambda_{\text{ex}})=0$
if $\theta=\frac{1}{1+\delta}\in (0,1)$.
\end{proof}
 
\begin{figure}[t]
  \centering
  \includegraphics[width=0.48\textwidth,clip]{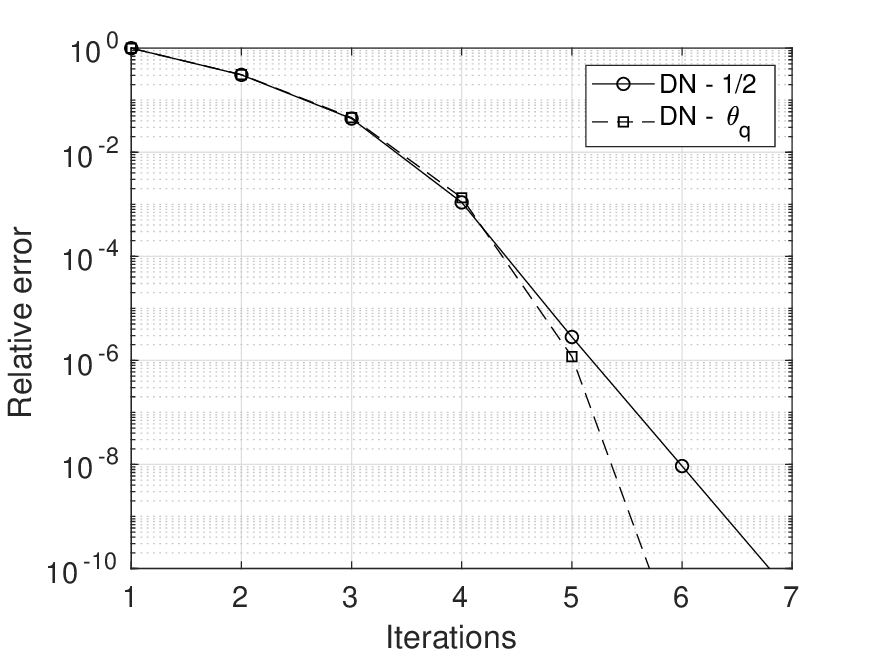}\quad
  \includegraphics[width=0.48\textwidth,clip]{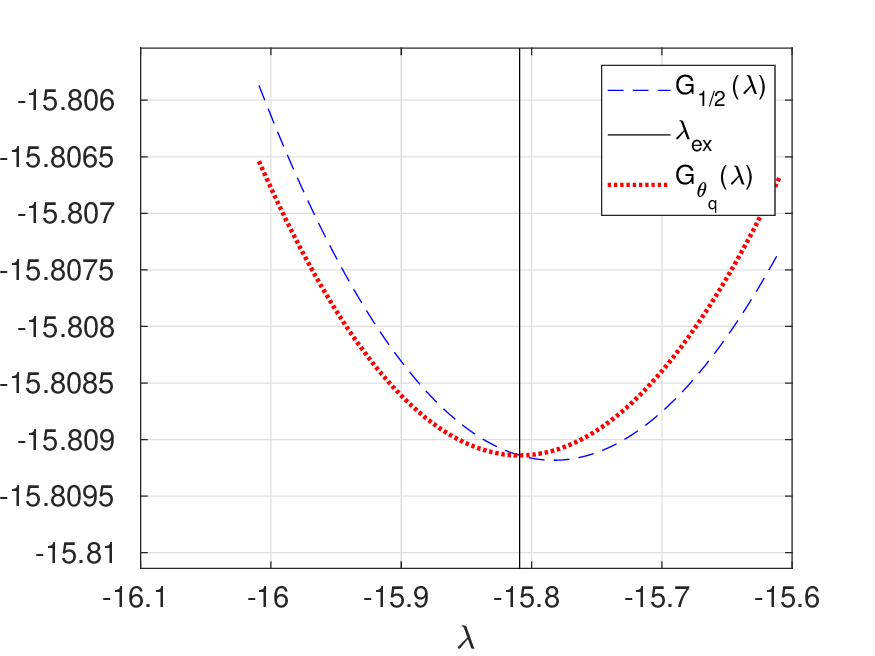}
  \includegraphics[width=0.48\textwidth,clip]{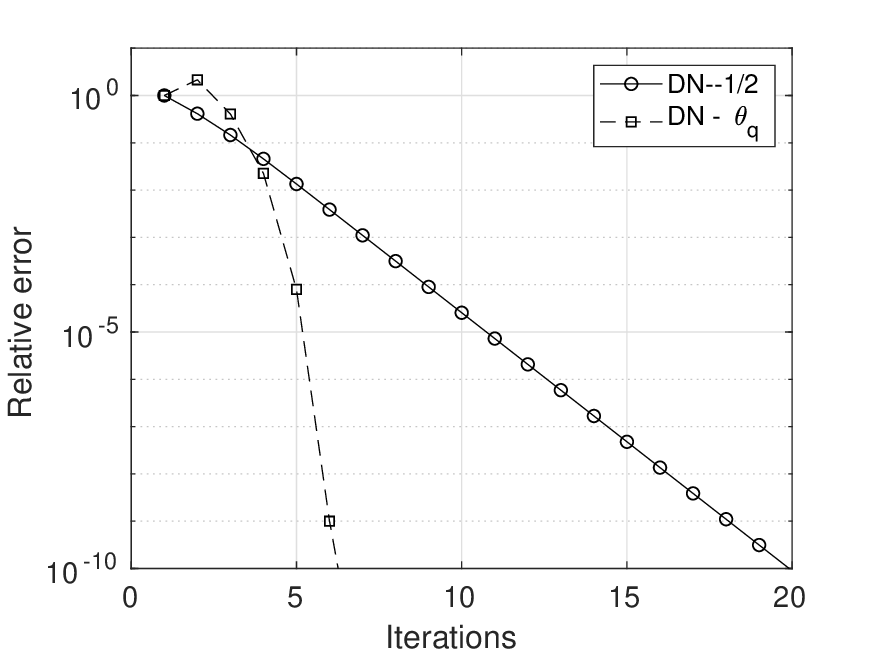}\quad
  \includegraphics[width=0.48\textwidth,clip]{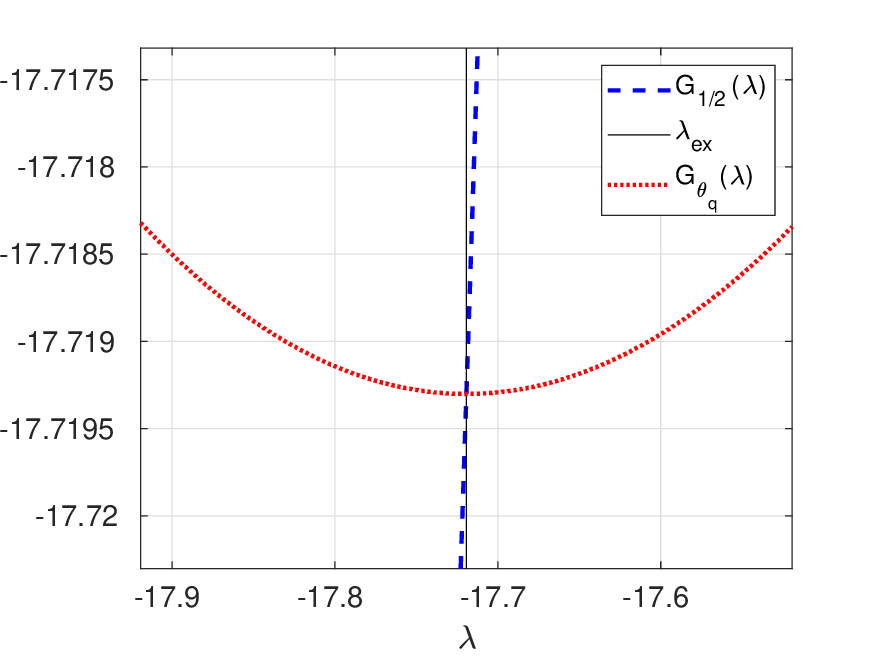}
  \caption{In the left panels, we show the convergence curves, and
  in the right panels we plot $G_{\theta}(\lambda)$. The top-row
  refers to a symmetric decomposition, and the bottom-row to an
  asymmetric one.}\label{kumbhar_p_mini_17_Fig:quadratic_convergence}
\end{figure}

To illustrate Theorem \ref{kumbhar_p_mini_17_thm:quadratic_convergence} numerically,
we consider $\mathcal{L}(u)=-\partial_x ((1+\alpha u^2)\partial_x u)$,
$\Omega=(0,1)$, $f(x)=100 x$, $u(0)=0$ and $u(1)=-20$. In the top-row
of Fig. \ref{kumbhar_p_mini_17_Fig:quadratic_convergence},
we set the interface $\Gamma$ to $x=1/2$. In the left panel,
we plot the convergence curves for $\theta=1/2$ and for
$\theta_{\text{q}}:=\frac{1}{1+\delta}$. In this setting, $\delta=1.006$ and
$\theta_{\text{q}}=0.498$, so due to the symmetry of the decomposition,
$\theta_{\text{q}}$ is still very close to $1/2$.  In the right panel, we
plot $G_\theta(\lambda)$ and see that as $\theta$ changes, the
minimum of $G_\theta(\lambda)$ moves, such that it is attained at
$\lambda=\lambda_{\text{ex}}$ for $\theta=\theta_{\text{q}}$.

Next, in the bottom row of Fig \ref{kumbhar_p_mini_17_Fig:quadratic_convergence}, we
consider the same equation and boundary conditions, but $\Gamma$ is
now at $x=0.3$. The decomposition is asymmetric, with $\delta=0.43$
and $\theta_{\text{q}}=0.699$.  The left panel shows clearly that for
$\theta=1/2$ the convergence is linear, while for $\theta=\theta_{\text{q}}$,
the DN method converges quadratically. In the right panel, we observe
that $G_{1/2}(\lambda)$ does not have a local extremum at
$\lambda=\lambda_{\text{ex}}$, while $G_{\theta_q}(\lambda)$ does.
Theorem 1 does not easily generalize to higher dimensions,
  since $\text{DtN}^\prime_j$ are then matrices, and 
  the relaxation parameter would have to be an
    operator. Numerically we observed for symmetric
  decompositions fast convergence for $\theta=0.5$, while for
  asymmetric decompositions, $\theta$ needs to be tuned for good
  performance.

\section{Mesh independent convergence}\label{kumbhar_p_mini_17_sec:mesh_ind}

\begin{figure}[t] 
  \centering
  \includegraphics[width=0.48\textwidth,clip]{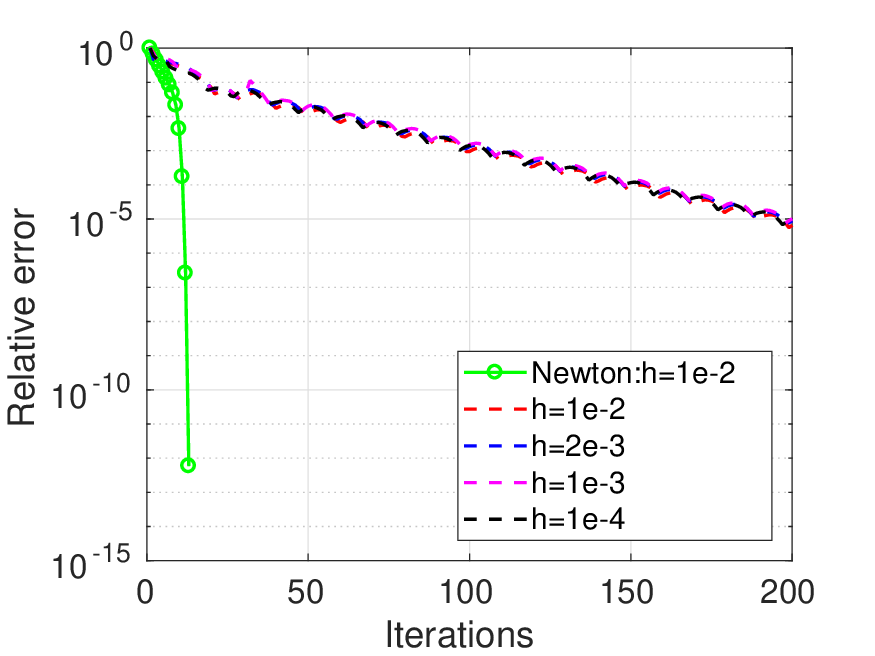}\quad
  \includegraphics[width=0.48\textwidth,clip]{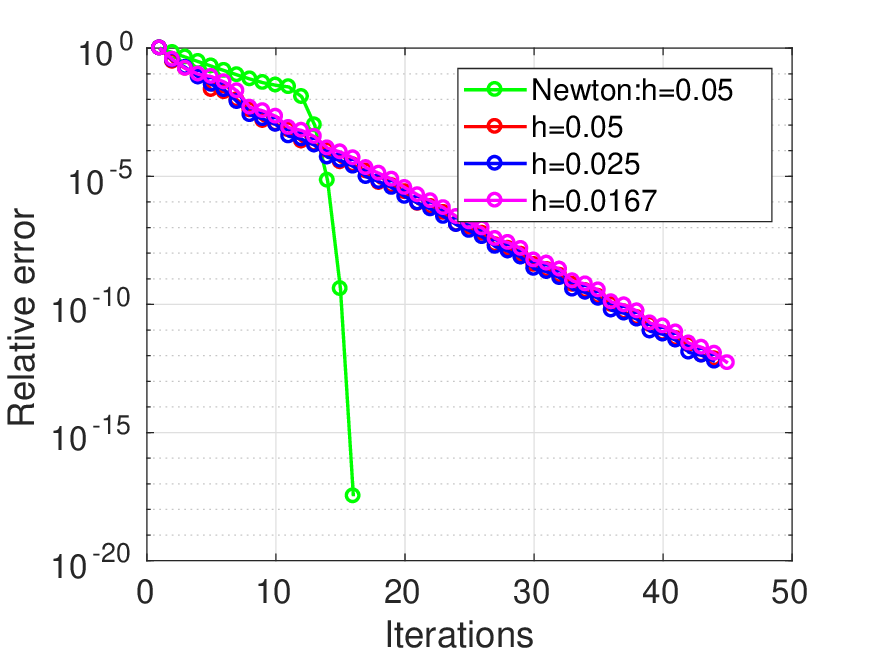}
  \caption{Convergence behavior of nonlinear DN for
    different mesh sizes in 1D (left) and 2D
    (right).}\label{kumbhar_p_mini_17_Fig:mesh_ind}
\end{figure}

One of the attractive features of the DN method for linear problems is
that it achieves mesh independent convergence. Does this also
  hold for the nonlinear DN method~\eqref{kumbhar_p_mini_17_eq:DN_method}?  We first
define the nonlinear DN method for multiple subdomains. Motivated
by the definition of the DN method for the linear case in
\cite{kumbhar_p_mini_17_scalabilty_CCGT}, we divide the domain $\Omega:=(0,L)\times
(0,L)$ into $N$ nonoverlapping subdomains
$\Omega_j=(\Gamma_{j-1},\Gamma_j)\times(0,L)$, with $\Gamma_0=0$ and
$\Gamma_N=L$.  The nonlinear DN method for multiple subdomains is then
defined for the interior subdomains by
\begin{equation*}
  \begin{array}{rcll}
    \mathcal{L}(u^n_j)&=&f_j \quad &\text{in $\Omega_j$}, \\
    \mathcal{N}_j u^n_j(\Gamma_{j-1},\cdot)&=&
    -\mathcal{N}_{j-1} u^{n}_{j-1}(\Gamma_{j-1},\cdot)  \quad
    &\text{on$\quad \Gamma_{j-1}$}, \\
    u^n_j(\Gamma_j)&=&(1-\theta)u^{n-1}_{j}(\Gamma_j,\cdot)
    + \theta u^{n-1}_{j+1}(\Gamma_j,\cdot)  \quad &\text{on $\Gamma_j$}, 
  \end{array}
\end{equation*}
where $\theta\in(0,1)$, and for the left and right most subdomains by
\begin{equation*}
\begin{array}{rclrcl}
  \mathcal{L}(u_1^n)&=&f_1,\quad \text{in}\quad \Omega_1, &
  \mathcal{L}(u_N^n)&=&f_N,\quad \text{in}\quad \Omega_N, \\
  u^n_1(\Gamma,\cdot)&=&g(0),&
  \hspace{-1em}\mathcal{N}_N u^n_N(\Gamma_{N-1},\cdot)&=&-\mathcal{N}_{N-1}
  u^{n}_{N-1}(\Gamma_{N-1},\cdot),\\
  u^n_1(\Gamma_1,\cdot)&=&(\!1\!-\!\theta\!)u^{n-1}_{1}(\Gamma_1,\cdot)
  \!+\! \theta  u^{n-1}_{2}(\Gamma_1,\cdot), & u^n_N(L,\cdot)&=&g(L).
\end{array}
\end{equation*}
We perform two experiments, one in 1D and one in 2D. For the
1D case, we consider the nonlinear diffusion equation
$-\partial_x\left((1+u^2)\partial_xu\right)=0$, with $u(0)=0$ and
$u(1)=20$. We divide the domain $\Omega=(0,1)$ into ten equal
subdomains. We then plot the relative error of the nonlinear DN for
four different mesh sizes $h$=1e-2, $h$=2e-3, $h$=1e-3, and
$h$=1e-4. The left plot in Fig. \ref{kumbhar_p_mini_17_Fig:mesh_ind}
shows that the convergence rate of the nonlinear DN is
independent of mesh size, while it is quadratic for Newton's Method. We repeat a similar experiment in 2D, but
now the domain $\Omega=(0,1)\times(0,1)$ is divided into four equal
subdomains. Even in 2D, we observe the mesh independent convergence of
the nonlinear DN method, see the right plot of Fig. \ref{kumbhar_p_mini_17_Fig:mesh_ind}.

\section{Dirichlet-Neumann Preconditioned Exact Newton (DNPEN)}

In Section \ref{kumbhar_p_mini_17_sec:Nil_quad}, we observed that under some special
conditions on the exact solution of the nonlinear problem and
$\theta$, the nonlinear DN method \eqref{kumbhar_p_mini_17_eq:DN_method} can be
  nilpotent. Moreover, the nonlinear DN method can also converge
quadratically.  But to achieve this, we need to tune the parameter
$\theta$ according to some a priori knowledge of the exact solution of
the nonlinear problem. Thus in general, the nonlinear DN method
converges linearly (as shown in Fig \ref{kumbhar_p_mini_17_Fig:mesh_ind}).
     
Iterative methods can be used as preconditioners to achieve
faster convergence, see \cite{kumbhar_p_mini_17_bookCG} for the linear
  case, and \cite{kumbhar_p_mini_17_gander2017origins} for a historical introduction
  including also the nonlinear case. It was proposed in
\cite{kumbhar_p_mini_17_dolean2016nonlinear, kumbhar_p_mini_17_SUBRAS} to use the nonlinear Restricted
Additive Schwarz (RAS) and nonlinear Substructured RAS (SRAS)
methods as preconditioner for Newton's method. We use the
same idea here and apply Newton's method to the fixed point
equation of the nonlinear DN method \eqref{kumbhar_p_mini_17_eq:DN_method_SUB},
which represents a systematic way of constructing non-linear
  preconditioners \cite{kumbhar_p_mini_17_gander2017origins}.  The fixed point version
of \eqref{kumbhar_p_mini_17_eq:DN_method_SUB} can be written as
\begin{equation}\label{kumbhar_p_mini_17_define:F_DN} 
\mathcal{F}(\mathbf{\lambda}):=\lambda-G(\lambda) =\theta\lambda - \theta \text{NtD}_2 \left(-\text{DtN}_1\left(\lambda,\psi_1\right),\psi_2\right).
\end{equation} 
Applying Newton to \eqref{kumbhar_p_mini_17_define:F_DN} we obtain a new method called
Dirichlet Neumann Preconditioned Exact Newton (DNPEN) method.
     
We saw in Section~\ref{kumbhar_p_mini_17_sec:Nil_quad} that the DN method can be
  nilpotent in certain cases. Can DNPEN still be nilpotent?  Let
$\lambda_{\text{ex}}$ denote the fixed point of the
iteration \eqref{kumbhar_p_mini_17_eq:DN_method_SUB}. Let us assume that the Dirichlet
Neumann method converges in one iteration. This means that $G$ defined
in~\eqref{kumbhar_p_mini_17_eq:DN_method_SUB} satisfies $\lambda_{\text{ex}}=G(\lambda^0)$
for any initial guess $\lambda^0$. This shows that the map $G$ is
constant, and hence $\mathcal{F}^\prime(\lambda)$ reduces to the
identity matrix. Moreover, one step of Newton's method applied
to~\eqref{kumbhar_p_mini_17_eq:DN_method_SUB} can then be written as
\[
     \lambda^1=\lambda^0-(\mathcal{F}^\prime(\lambda^0))^{-1}\mathcal{F}(\lambda^0)=\lambda^0-\mathcal{F}(\lambda^0) =G(\lambda^0)=\lambda_{\text{ex}},
\]
and hence DNPEN will also be nilpotent in that case. We further have also the following result.
\begin{theorem}\label{kumbhar_p_mini_17_th:DNPEN}
  The convergence of DNPEN does not depend on the relaxation parameter
$\theta$ in the DN preconditioner.
\end{theorem}     
\begin{proof}
  The function $\mathcal{F}$ from \eqref{kumbhar_p_mini_17_define:F_DN} corresponding to
  DNPEN can we rewritten as $\mathcal{F}(\lambda)=\theta
  \mathcal{K}(\lambda,\psi_1,\psi_2)$, where
  $\mathcal{K}(\lambda,\psi_1,\psi_2):=\lambda - \text{NtD}_2
  \left(-\text{DtN}_1\left(\lambda,\psi_1\right),\psi_2\right)$. Thus,
  Newton's iteration reads
  \[\lambda^{k+1}\!=\!\lambda^k - \left(J\mathcal{F}(\lambda^k)\right)^{-1}
  \!\!\mathcal{F}(\lambda^k)\!=\!\lambda^k- \left(\theta J\mathcal{K} (\lambda^k)\right)^{-1}\!\!\theta \mathcal{K}(\lambda^k)\!=\!\lambda^k- \left(J\mathcal{K} (\lambda^k)\right)^{-1}\!\!\mathcal{K}(\lambda^k),
  \]
  which shows that the Newton correction does not depend on the
  relaxation parameter $\theta$. The iterates of Newton's method will
  thus only depend on $\mathcal{K}$, and DNPEN has $\theta$
  independent convergence.
\end{proof}	

The above theorem shows that when using DNPEN, one does not need
to search for an optimal choice of $\theta$, in contrast to the
nonlinear DN method \eqref{kumbhar_p_mini_17_eq:DN_method}.
        
\begin{figure}[t]
\centering
\includegraphics[width=0.48\textwidth,clip]{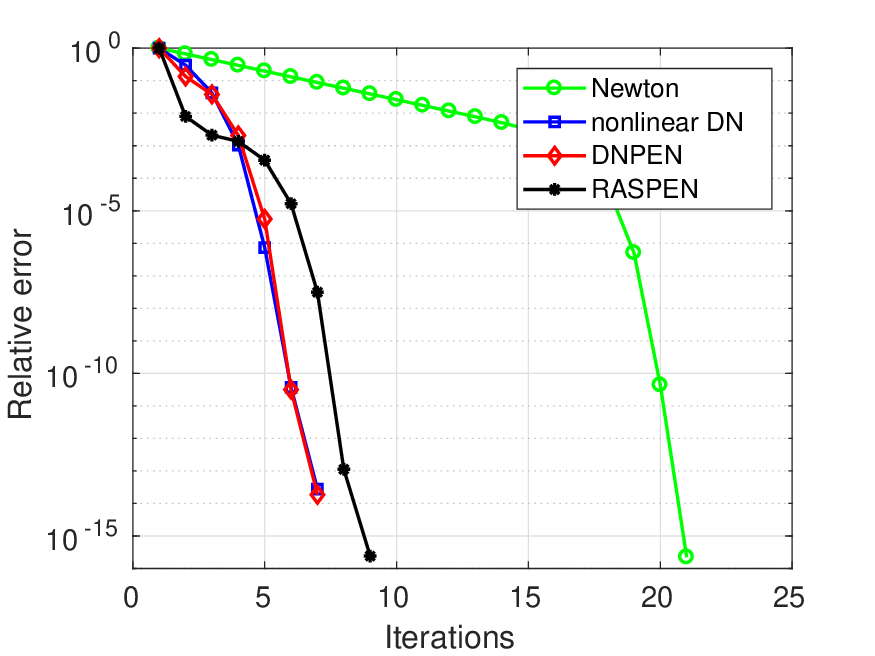}\quad
\includegraphics[width=0.48\textwidth,clip]{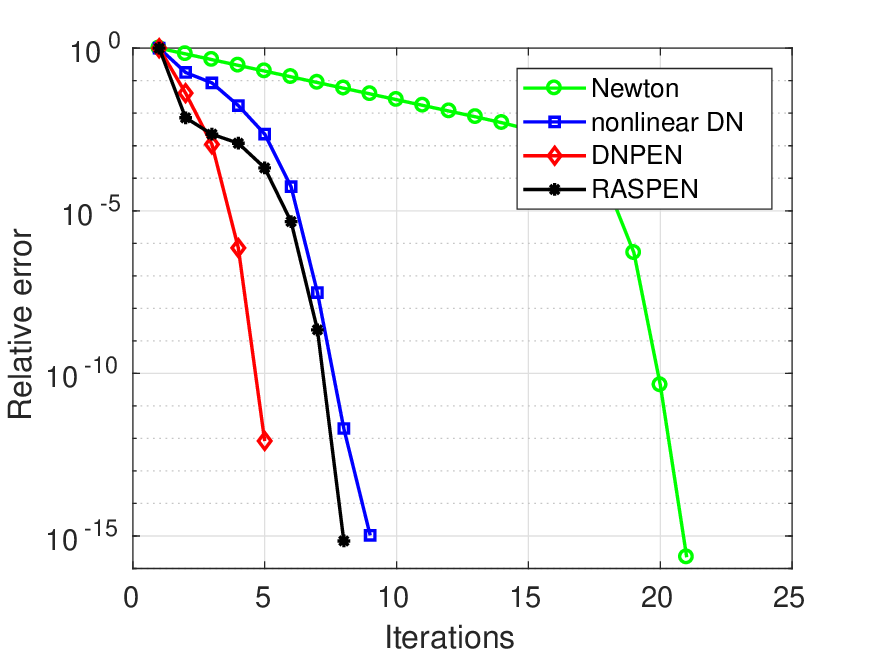}
\caption{Comparison of DNPEN (with optimal $\theta$) with
  unpreconditioned Newton, nonlinear DN (with optimal $\theta$)
  and RASPEN for a symmetric partition (left) and an asymmetric
  partition (right).}\label{kumbhar_p_mini_17_Fig:DNPEN_optimal}
\end{figure}

We now compare the convergence of DNPEN, the unpreconditioned
  Newton method, the nonlinear DN method \eqref{kumbhar_p_mini_17_eq:DN_method} and
RASPEN \cite{kumbhar_p_mini_17_dolean2016nonlinear}.
We consider the nonlinear diffusion problem $-\partial_x\left(\left(1+u^2\right)\partial_x u\right)=f$
on $\Omega=(0,1)$ decomposed into two equally sized
subdomains, with $u(0)=0$, $u(1)=10$ and $f(x)=\sin(10\pi x)$.  For
both DN and DNPEN, we choose the optimal relaxation parameter provided
in Theorem~\ref{kumbhar_p_mini_17_thm:quadratic_convergence}.  The left plot in
Fig. \ref{kumbhar_p_mini_17_Fig:DNPEN_optimal}
shows that the iterative DN converges quadratically using the
optimal parameter and is very similar to DNPEN with no significant
gain in the number of iterations. The convergence curves also show
that the unpreconditioned Newton method is slower than all
preconditioned ones, and DNPEN has a slight advantage over
RASPEN.

We repeat the same experiment but now using an asymmetric partition of
the domain $\Omega$. The right plot in
Fig. \ref{kumbhar_p_mini_17_Fig:DNPEN_optimal} shows that for this configuration, DNPEN
is the fastest while again unpreconditioned Newton is the slowest
among the methods considered. Moreover, DNPEN is significantly
faster than the nonlinear DN method.

\begin{figure}
\centering
\includegraphics[width=0.4625\textwidth,clip]{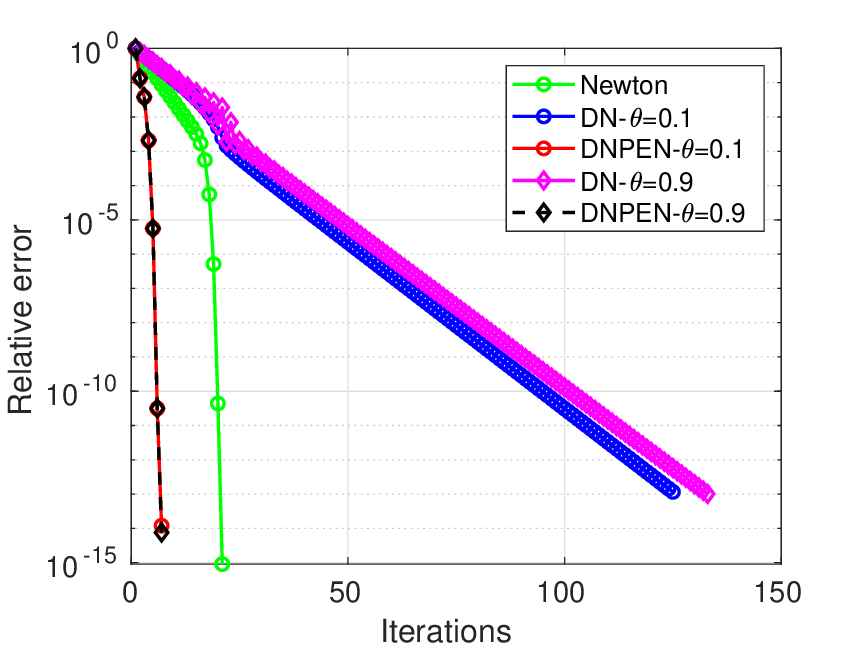}\quad
\includegraphics[width=0.4625\textwidth,clip]{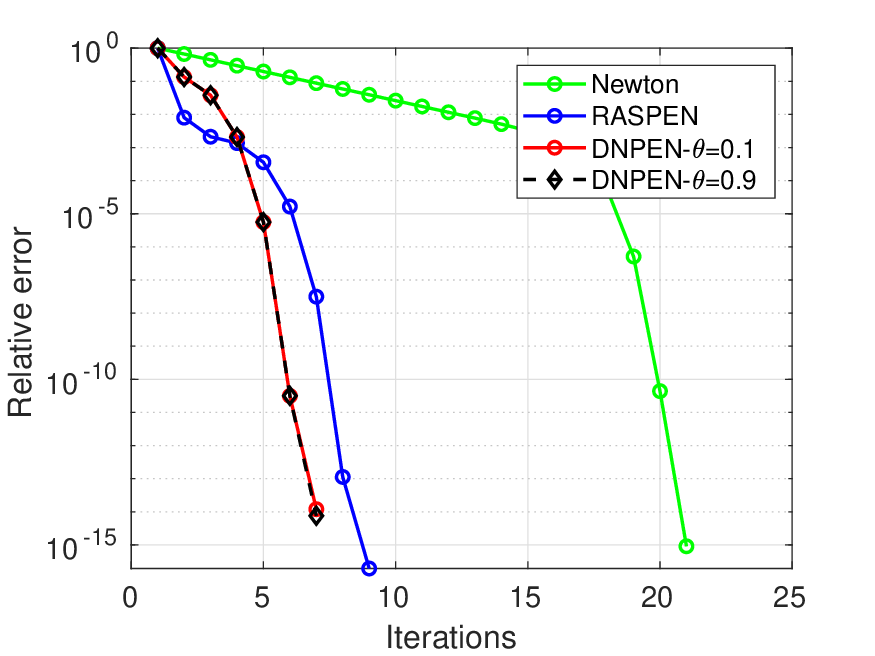}
\caption{Comparison of DNPEN with the unpreconditioned Newton
  method and nonlinear DN (left) and with RASPEN (right) for two
  different non optimal $\theta$'s.}\label{kumbhar_p_mini_17_Fig:DNPEN_non_optimal}
\end{figure}

Finally, we illustrate numerically that the convergence of DNPEN
does not depend on $\theta$. We know that in general, the nonlinear DN
method converges linearly, and it is not always possible to find an
optimal $\theta$ such that it converge quadratically. We again
consider the symmetric partition of the domain and use the same
boundary conditions and force term as above. However, instead of
the optimal $\theta$, we consider two non-optimal $\theta$'s,
namely $\theta=0.1$ and $\theta=0.9$. The left plot in Fig.
\ref{kumbhar_p_mini_17_Fig:DNPEN_non_optimal}
shows the linear convergence of nonlinear DN for both $\theta=0.1$,
and $\theta=0.9$, and both are slower than the unpreconditioned
Newton method. However, DNPEN converges much faster than Newton's
method and in the same number of iterations for the two
different values $\theta=0.1$ and $\theta=0.9$. The right plot in
Fig. \ref{kumbhar_p_mini_17_Fig:DNPEN_non_optimal} shows that DNPEN is still faster
than RASPEN for both values $\theta$ considered.
     
\section{Conclusion}     

While iterative DN methods are known to converge linearly, we
proved that one can obtain
quadratic converge for some one-dimensional nonlinear problems and for a well chosen relaxation parameter
  $\theta$. Under specific conditions, the nonlinear DN method can
also become a direct solver, like in the linear case. We then extended
DN to multiple subdomains and numerically showed that its convergence
is mesh independent. We finally introduced the nonlinear
preconditioner DNPEN, proved that the convergence of DNPEN does not
depend on the relaxation parameter $\theta$, and observed
numerically that DNPEN is faster than unpreconditioned Newton,
nonlinear DN and RASPEN in all our examples.

\section*{Acknowledgements}
The third author acknowledges financial support from the Deutsche Forschungsgemeinschaft (DFG, German Research Foundation) - Project-ID 258734477 - SFB 1173. 

\bibliographystyle{plain}
\bibliography{Kumbhar_p_mini_17.bib}
\end{document}